\documentclass[a4paper,10pt]{article}
\usepackage{amsmath}
\usepackage{amssymb}
\usepackage{graphics}
\usepackage{rotating}
\usepackage{cite}
\usepackage{color}

\title{A generalization of Fueter's theorem}
\author{Daniel Alay\'{o}n-Solarz \underline{(danieldaniel@gmail.com)}}

\begin{document}

\maketitle

\begin{abstract}
We show that Fueter's theorem holds for a more general class of quaternionic functions than those constructed by the Fueter's method.\end{abstract}
\section*{Description of result}

Let $f(z)$ be a holomorfic function defined on the upper complex plane. If we write $z=x+iy$ and $f(z)=u(z)+iv(z)$ we can construct a quaternionic function $f(p)=u(t,r)+v(t,r)\iota$ where $p=t+xi+yj+zk$, $r=\sqrt{x^2+y^2+z^2}$ and

\begin{equation}
\iota=\frac{xi+yj+zk}{\sqrt{x^2+y^2+z^2}}
\end{equation}

If we denote $D_{l}$ and $D_{r}$ the left and right Fueter (or Dirac) operator then Fueter's theorem asserts that

 \begin{equation}
 D_{l}\Delta{f}=D_{r}\Delta{f}=0
\end{equation}

such functions are never regular unless they are real constants, instead they satisfy:

 \begin{equation}
 D_{l}(f)=D_{r}(f)= \frac{-2v(t,r)}{r}
\end{equation}

One can easily prove Fueter's theorem observing that:

\begin{equation}
- \frac{1}{2}\Delta(f)= \bar{D}(\frac{v}{r})= -\frac{\iota}{r}\frac{\partial f}{\partial t}+\frac{\iota}{r^2}v(t,r)
 \end{equation}
 
 which is an axial symmetric function that can be written in the form $\tilde{f}=\tilde{u}+\iota\tilde{v}$ and that such functions are regular if and only if they satify
 
 \begin{equation}
(\frac{\partial}{\partial t}+\iota \frac{\partial}{\partial r})\tilde{f}=\frac{2\tilde{v}}{r}
  \end{equation}
  
  Our result points out that if a quaternionic function of Class $C^{2}$ can be written as $f(p)=u(p)+\iota v(p)$ and satisfy
  
  \begin{equation}
 D_{l}(f)=\frac{-2v}{r}
\end{equation}
  
  then its laplacian will be left- and right-regular. Note that we are dropping the condition that the function is obtained by the Fueter's method and so the resulting functions need not be axial symmetric.
  
%%%%%%%%%%%%%%%%%%%%%%%%%%%%%%%%%%%%%%%%%%%%%%%%%%%%%%%%%%%%%%%%%%%%%%%%%%%%%%%
%                                SECTION II
%%%%%%%%%%%%%%%%%%%%%%%%%%%%%%%%%%%%%%%%%%%%%%%%%%%%%%%%%%%%%%%%%%%%%%%%%%%%%%%
\section*{II. Preliminaries}
Let $f(p)$ be a quaternionic, $C^{2}$ function that satisfy $f(p)p = pf(p)$. Then there exists real functions $u,v$ such that 

\begin{equation}
f(p) = u(p) + \iota v(p)
\end{equation} 
Notice that complex functions obtained by the Fueter's method are of this form. We write the quaternion $p$ as
\begin{equation}
p = t + r \iota
\end{equation}
And we parametrize  $\iota$ by spherical coordinates.
\begin{equation}
\iota = (\cos\alpha \sin\beta, \sin\alpha \sin\beta, \cos \beta)
\end{equation}

Let $D_{]}$ denote the left-Fueter operator, so
\begin{equation}
D_{l}:= \frac{\partial }{\partial t}+ i\frac{\partial }{\partial x} + j\frac{\partial }{\partial y} + k\frac{\partial }{\partial z} 
\end{equation}
And we are interested in functions of the form (7) that satisfy
 \begin{equation}
 D_{l}f = \frac{-2v}{r}
 \end{equation}
The Fueter operator in $(t,r,\alpha,\beta)$ coordinates is written as
\begin{equation}
D_{l} = \frac{\partial }{\partial t} + \iota \frac{\partial }{\partial r}  - \frac{1}{r} \frac{\partial}{\partial_{l} \iota}
\end{equation}

where the symbol divided by $r$ is defined as:

\begin{equation}
\frac{\partial}{\partial_{l} \iota} := {\iota}_{\alpha}^{-1}\frac{\partial}{\partial \alpha} +  {\iota}_{\beta}^{-1}\frac{\partial}{\partial \beta}
\end{equation}

and $\iota_{\alpha}$ and $\iota_{\beta}$ denote the partial derivative of $\iota$ respect of $\alpha$,$\beta$

As a function that satisfy (11) is holomorphic on the variables $t,r$ the we conclude that it must hold

 \begin{equation}
 \frac{\partial f}{\partial_{l} \iota} = 2v
 \end{equation}
 
 which occurs if and only if $u$ and $v$ satisfy the following Cauchy-Riemann type system:
  \begin{equation}
 \frac{1}{\sin \beta}\frac{\partial u}{\partial \alpha} =  \frac{\partial v}{\partial \beta}
 \end{equation}
 \begin{equation}
 \frac{1}{\sin \beta}\frac{\partial v}{\partial \alpha} =  -\frac{\partial u}{\partial \beta}
 \end{equation}
%%%%%%%%%%%%%%%%%%%%%%%%%%%%%%%%%%%%%%%%%%%%%%%%%%%%%%%%%%%%%%%%%%%%%%%%%%%%%%%
%                                SECTION III
%%%%%%%%%%%%%%%%%%%%%%%%%%%%%%%%%%%%%%%%%%%%%%%%%%%%%%%%%%%%%%%%%%%%%%%%%%%%%%%
\section*{II. Proof of result}

We calculate explicitly that if $f= u + \iota v$ satisfy (11) then:

\begin{equation}
\nabla \frac{v}{r}=0
\end{equation}

Observe that
\begin{equation}
 \bar{D_{l}} (\frac{v}{r}) = -\frac{\iota}{r}\frac{\partial f}{\partial t}+\frac{\iota}{r^2}v + \frac{1}{r^2} \frac{\partial v}{\partial_{l} \iota}
 \end{equation} 
 
 Before proceeding, we will remark some properties that will be helpful.
 
 The first remark is that if $f$ satisfy (11) then $\frac{\partial{f}}{\partial t}$ will also satify (11). The second remark is to note that, for all functions $f$
 \begin{equation}
 \frac{\partial }{\partial_{l} \iota}(\iota f) = 2f -\iota \frac{\partial f }{\partial_{l} \iota}
 \end{equation}
 The third remark is that if $f= u + \iota v$
 \begin{equation}
 \frac{\partial f}{\partial_{l} \iota} = 2v 
 \end{equation}
 is equivalent to
 \begin{equation}
 \frac{\partial u}{\partial_{l} \iota} = \iota \frac{\partial v}{\partial_{l} \iota}
 \end{equation}
 The fourth remark is that if
 \begin{equation}
 \frac{\partial f}{\partial_{l} \iota} = 2v 
 \end{equation}
 then 
 \begin{equation}
 \frac{\partial }{\partial_{l} \iota}(\iota f) = 2u 
 \end{equation}
 
 Now, for the sake of readability we will apply  $D_{l}$ only to the first summand in (18) and we have:
 \begin{equation}
 D_{l}(-\frac{\iota}{r}\frac{\partial f}{\partial t}) = \frac{1}{r^2}(-\frac{\partial f}{\partial t} + \frac{\partial }{\partial_{l} \iota}(\iota \frac{\partial f}{\partial t})) = \frac{1}{r^2}(-\frac{\partial f}{\partial t} + 2\frac{\partial u}{\partial t})
 \end{equation}
 \begin{equation}
 = \frac{1}{r^2}(\frac{\partial u}{\partial t}- \iota \frac{\partial v}{\partial t})
 \end{equation}
 We apply  $D_{l}$ to the second summand in (18) and we have:
 \begin{equation}
  D_{l}(\frac{\iota}{r^2}v) = \frac{\iota}{r^2}\frac{\partial v}{\partial t} + \iota (-2\frac{ \iota v}{r^3}+ \frac{\iota}{r^2}\frac{\partial v}{\partial r}) + \frac{1}{r^3}\frac{\partial }{\partial_{l}\iota}(\iota v)
 \end{equation}
 \begin{equation}
 =\frac{1}{r^2}(\iota \frac{\partial v}{\partial t} - \frac{\partial v}{\partial r}) + \frac{1}{r^3}\iota \frac{\partial v}{\partial_{l}\iota}
 \end{equation}
 We apply $D_{l}$ to the third summand in (18) and we have:
 \begin{equation}
  D_{l}(\frac{1}{r^2}\frac{\partial v}{\partial_{l}\iota}) = \frac{1}{r^2}\frac{\partial}{\partial t}\frac{\partial v}{\partial_l \iota} + \iota (\frac{-2}{r^3}\frac{\partial v}{\partial_{l} \iota}+\frac{1}{r^2} \frac{\partial }{\partial r} \frac{\partial v}{\partial_{l} \iota})- \frac{1}{r^3}\frac{\partial^{2} v}{\partial_{l} \iota^{2}} =
 \end{equation}
  \begin{equation}
  \frac{1}{r^2}(\frac{\partial}{\partial t}\frac{\partial v}{\partial_{l} \iota}+\iota \frac{\partial}{\partial r}\frac{\partial v}{\partial_{l}\iota}) + \frac{1}{r^3}(-2\iota\frac{\partial v}{\partial_{l}\iota}-\frac{\partial^2 v}{\partial_{l}\iota^2})
  \end{equation}
  
  Now observe that in the first summand in (29)
  \begin{equation}
  \frac{\partial}{\partial t}\frac{\partial v}{\partial_{l} \iota}+\iota \frac{\partial}{\partial r}\frac{\partial v}{\partial_{l}\iota} = \frac{\partial}{\partial t}(\frac{\partial v}{\partial_{l} \iota}+ \iota \frac{\partial u}{ \partial_{l} \iota}) = 0
  \end{equation}
 because of the third remark. We also observe that the sum of the first two summands in (25) and (27) is zero because $f$ is holomorphic in $t,r$. So, for now we have established that 
 \begin{equation}
 \nabla \frac{v}{r} = \frac{1}{r^3}(\iota \frac{\partial v}{\partial_{l}\iota}  -2\iota\frac{\partial v}{\partial_{l}\iota}-\frac{\partial^2 v}{\partial_{l}\iota}) = \frac{1}{r^3}(-\iota\frac{\partial v}{\partial_{l}\iota}-\frac{\partial^2 v}{\partial_{l}\iota^2})
 \end{equation}
 so we now must show that
 \begin{equation}
 \frac{\partial^2 v}{\partial_{l}\iota}=-\iota\frac{\partial v}{\partial_{l}\iota}
 \end{equation}
 We first apply the definition
 \begin{equation}
 \frac{\partial^2 v}{\partial_{l}\iota} = ( {\iota}_{\alpha}^{-1}\frac{\partial}{\partial \alpha} +  {\iota}_{\beta}^{-1}\frac{\partial}{\partial \beta})( {\iota}_{\alpha}^{-1}\frac{\partial v}{\partial \alpha} +  {\iota}_{\beta}^{-1}\frac{\partial v}{\partial \beta})
 \end{equation}
 to obtain
 \begin{equation}
 \iota^{-1}_{\alpha}((\iota^{-1}_{\alpha})_{\alpha}\frac{\partial v}{\partial \alpha}+\iota^{-1}_{\alpha}\frac{\partial^{2} v}{\partial \alpha^{2}}+(\iota^{-1}_{\beta})_{\alpha}\frac{\partial v}{\partial \beta}+\iota^{-1}_{\beta}\frac{\partial^2 v}{\partial \alpha \partial \beta}) +
 \end{equation}
 \begin{equation}
 \iota^{-1}_{\beta}((\iota^{-1}_{\alpha})_{\beta}\frac{\partial v}{\partial \alpha}+\iota^{-1}_{\alpha}\frac{\partial^2 v}{\partial \beta \partial \alpha}+(\iota^{-1}_{\beta})_{\beta}\frac{\partial v}{\partial \beta}+\iota^{-1}_{\beta}\frac{\partial^2 v}{\partial \beta^2})
 \end{equation}
 the cross derivatives vanish because $v$ is assumed of class $C^{2}$ and $\iota_{\alpha}$ and $\iota_{\beta}$ anticommute, so we have left:
 \begin{equation}
  \iota^{-1}_{\alpha}((\iota^{-1}_{\alpha})_{\alpha}\frac{\partial v}{\partial \alpha}+\iota^{-1}_{\alpha}\frac{\partial^{2} v}{\partial \alpha^{2}}+(\iota^{-1}_{\beta})_{\alpha}\frac{\partial v}{\partial \beta}) + 
 \end{equation}
 \begin{equation}
 \iota^{-1}_{\beta}((\iota^{-1}_{\alpha})_{\beta}\frac{\partial v}{\partial \alpha}+(\iota^{-1}_{\beta})_{\beta}\frac{\partial v}{\partial \beta}+\iota^{-1}_{\beta}\frac{\partial^2 v}{\partial \beta^2}
 \end{equation}
 However, we have the following equalities, which can be easily verified
 \begin{equation}
  \iota^{-1}_{\alpha}(\iota^{-1}_{\alpha})_{\alpha}+\iota^{-1}_{\beta}(\iota^{-1}_{\alpha})_{\beta} = -\iota \iota^{-1}_{\alpha}
 \end{equation}
 \begin{equation}
 \iota^{-1}_{\beta}(\iota^{-1}_{\beta})_{\beta} = -\iota \iota^{-1}_{\beta}
 \end{equation} 
 \begin{equation}
 \iota^{-2}_{\beta} = -1
 \end{equation}
 \begin{equation}
 \iota^{-2}_{\alpha} = -\frac{1}{\sin^{2} \beta}
 \end{equation}
 \begin{equation}
 \iota^{-1}_{\alpha}(\iota^{-1}_{\beta})_{\alpha} = -\cot \beta
 \end{equation}
 so we have that for any function $v \in C^{2}$ holds
 \begin{equation}
  \frac{\partial^2 v}{\partial_{l}\iota^{2}} = -\iota \frac{\partial v}{\partial_{l}\iota} - \frac{1}{\sin^2 \beta }\frac{\partial^2 v}{\partial \alpha^2} - \frac{\partial^2 v}{\partial \beta^2} - \cot \beta \frac{\partial v}{\partial \beta}
  \end{equation}
  So in order for our result to hold, we must have
  \begin{equation}
   \frac{1}{\sin^2 \beta }\frac{\partial^2 v}{\partial \alpha^2} +\frac{\partial^2 v}{\partial \beta^2} = -\cot \beta \frac{\partial v}{\partial \beta}
  \end{equation}
  but this can be easily deduced from (15) and (16) together.
  We have determined that if a function $f$ is such that $u,v$ are $C^{2}$ and satisfy (6) then its laplacian will be left regular, we now prove that its laplacian will be right-regular also. Fortunately, the hard word is already done and so this is almost immediate.
 First we observe that the operators $D_{l}$,$D_{r}$,$\bar{D_{l}}$ and $\bar{D_{r}}$ all commute with each other. It is also true that for any scalar function, say, $g$ we have:
\begin{equation}
D_{l}g=D_{r}g
\end{equation}
Now we use the fact that the right hand side of (6) is an scalar function. We already know that
\begin{equation}
0= D_{l}\bar{D_{l}}D_{l} f = D_{l}\bar{D_{l}}(-2\frac{v}{r}) = \bar{D_{l}}D_{l}(-2\frac{v}{r}) = \bar{D_{l}}D_{r}(-2\frac{v}{r}) 
\end{equation}
\begin{equation}
=D_{r}\bar{D_{l}}(-2\frac{v}{r}) = D_{r}\bar{D_{l}}D_{l}f
\end{equation}

And so have we proved our

\newtheorem{prop1}{Theorem}
\begin{prop1}
Let $f$ be a quaternionic, $C^{2}$ function that can be written as $f=u + \iota v$, where $u,v$ are real functions and that satisfy
\begin{equation}
D_{l}f = \frac{-2v}{r}
\end{equation} 
Then it holds that
\begin{equation}
 D_{l}\Delta{f}=D_{r}\Delta{f}=0
\end{equation}
\end{prop1}

\textbf{Remark} The author suspects that the condition $C^{2}$ can be relaxed to be $C^{1}$, and that the $C^{2}$ can be deduced from the holomorphism in $t,r$. But he doesn't have a proof yet.

%%%%%%%%%%%%%%%%%%%equation%%%%%%%%%%%%%%%%%%%%%%%%%%%%%%%%%%%%%%%%%%%%%%%%%%%%%%%%%%%%
%                                SECTION IV
%%%%%%%%%%%%%%%%%%%%%%%%%%%%%%%%%%%%%%%%%%%%%%%%%%%%%%%%%%%%%%%%%%%%%%%%%%%%%%%

%%%%%%%%%%%%%%%%%%%%%%%%%%%%%%%%%%%%%%%%%%%%%%%%%%%%%%%%%%%%%%%%%%%%%%%%%%%%%%%
%                                SECTION V
%%%%%%%%%%%%%%%%%%%%%%%%%%%%%%%%%%%%%%%%%%%%%%%%%%%%%%%%%%%%%%%%%%%%%%%%%%%%%%%

%%%%%%%%%%%%%%%%%%%%%%%%%%%%%%%%%%%%%%%%%%%%%%%%%%%%%%%%%%%%%%%%%%%%%%%%%%%%%%%
%                                SECTION VI
%%%%%%%%%%%%%%%%%%%%%%%%%%%%%%%%%%%%%%%%%%%%%%%%%%%%%%%%%%%%%%%%%%%%%%%%%%%%%%%

%%%%%%%%%%%%%%%%%%%%%%%%%%%%%%%%%%%%%%%%%%%%%%%%%%%%%%%%%%%%%%%%%%%%%%%%%%%%%%%
%                                SECTION VII
%%%%%%%%%%%%%%%%%%%%%%%%%%%%%%%%%%%%%%%%%%%%%%%%%%%%%%%%%%%%%%%%%%%%%%%%%%%%%%%

%%%%%%%%%%%%%%%%%%%%%%%%%%%%%%%%%%%%%%%%%%%%%%%%%%%%%%%%%%%%%%%%%%%%%%%%%%%%%%%
%                                REFERENCES
%%%%%%%%%%%%%%%%%%%%%%%%%%%%%%%%%%%%%%%%%%%%%%%%%%%%%%%%%%%%%%%%%%%%%%%%%%%%%%%

\end{document}